\documentclass[amscd,verbatim,12pt]{amsart}

\usepackage{tikz}
\usepackage{amssymb,amsmath,latexsym,graphics,youngtab,supertabular}
\usepackage{amsbsy}
\usepackage{amscd}
\usepackage{amsfonts}
\usepackage{amsthm}
\usepackage{epsfig}
\usepackage{times}
\usepackage{mathrsfs}
\usepackage{verbatim}
\usepackage{url}
\usepackage{txfonts}

\usepackage[english]{babel}
\usepackage[T1]{fontenc}
\usepackage[utf8]{inputenc} 

\usetikzlibrary{backgrounds,fit,decorations.pathreplacing}
\usetikzlibrary{arrows,shapes,positioning}
\usetikzlibrary{calc,through,backgrounds,chains}

\usetikzlibrary{arrows,shapes,snakes,automata,backgrounds,petri}

\usepackage[version=3]{mhchem}

\newtheorem{theorem}{Theorem}[section]

\theoremstyle{definition}
\newtheorem{definition}[theorem]{Definition}

\begin{document}

\title[Broken bracelets]{Broken bracelets, Molien series,
paraffin wax and an elliptic curve of conductor $48$}

\author[T. Amdeberhan]{Tewodros Amdeberhan}
\address{Department of Mathematics,
Tulane University, New Orleans, LA 70118}
\email{tamdeber@tulane.edu}

\author{Mah\.{i}r B\.{i}len Can}
\address{Department of Mathematics,
Tulane University, New Orleans, LA 70118}
\email{mcan@tulane.edu}

\author{Victor H. Moll}
\address{Department of Mathematics,
Tulane University, New Orleans, LA 70118}
\email{vhm@math.tulane.edu}


\date{\today}


\begin{abstract}
This paper introduces the concept of necklace binomial coefficients
motivated by the enumeration of a special type of sequences. Several properties 
of these coefficients are described, including a connection 
between their roots and an elliptic curve. Further links are given to a 
physical model from quantum mechanical supersymmetry as well as properties 
of alkane molecules in chemistry. 
\end{abstract}

\maketitle

\newcommand{\nn}{\nonumber}
\newcommand{\ba}{\begin{eqnarray}}
\newcommand{\ea}{\end{eqnarray}}
\newcommand{\realpart}{\mathop{\rm Re}\nolimits}
\newcommand{\imagpart}{\mathop{\rm Im}\nolimits}
\newcommand{\lcm}{\operatorname*{lcm}}
\newcommand{\vGam}{\varGamma}

\newtheorem{Definition}{\bf Definition}[section]
\newtheorem{Thm}[Definition]{\bf Theorem}
\newtheorem{Theorem}[Definition]{\bf Theorem}
\newtheorem{Example}[Definition]{\bf Example}
\newtheorem{Lem}[Definition]{\bf Lemma}
\newtheorem{Note}[Definition]{\bf Note}
\newtheorem{Cor}[Definition]{\bf Corollary}
\newtheorem{Prop}[Definition]{\bf Proposition}
\newtheorem{Conj}[Definition]{\bf Conjecture}
\newtheorem{Problem}[Definition]{\bf Problem}
\numberwithin{equation}{section}

\section{Introduction} \label{sec-intro}
\setcounter{equation}{0}

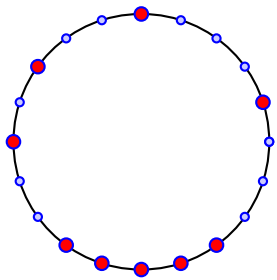
\begin{figure}[htp]

\begin{center}
    \begin{tikzpicture}[scale=1.7,cap=round,>=latex]
     
        \draw[thick] (0cm,0cm) circle(1cm);

    \foreach \x in {0,18,...,360} {
    \fill[draw=blue,fill=blue!20!,thick] (\x:1cm) circle(0.9pt);
    }
	\foreach \x in {234,252, 270,288,306 } {
    \fill[draw=blue,fill=red,thick] (\x:1cm) circle(1.5pt);
    }
	\foreach \x in {18, 90, 144,180 } {
    \fill[draw=blue,fill=red,thick] (\x:1cm) circle(1.5pt);
    }

    \end{tikzpicture}

\end{center}

\caption{A necklace with a medallion.}
\label{fig-med}
\end{figure}

A jeweler is asked to design a necklace consisting of a chain with $n$ 
placements for $k$ pieces of diamond.  The client ask for one group of 
$r$ diamonds to be placed next to each other and the remaining diamonds are 
to be isolated, that is, each one is mounted so that the two adjacent 
places are left empty. These special diamonds are  called the 
{\em medallion} of the necklace. Figure \ref{fig-med} shows a necklace of 
length $20$, with a medallion of length $5$ and four extra diamonds.

\begin{figure}[htp]
	
	\begin{center}
    \begin{tikzpicture}[scale=1,cap=round,>=latex]
    
     \path[draw,-,thick] (0,0) -- (9.5,0);
   
	\foreach \x in {0,1,...,19 } {
       \fill[draw=blue,fill=blue!20!,thick] (\x/2 ,0) circle(1.5pt);
    }

	\foreach \x in {0,1,5,9,12,14,17,18,19} {
       \fill[draw=blue,fill=red,thick] (\x/2 ,0) circle(2.5pt);
    }

    \end{tikzpicture}
	\end{center}
	\caption{A configuration.}
        \label{fig-1}
	\end{figure}
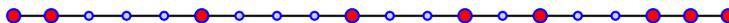

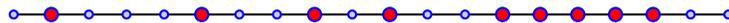
\begin{figure}[htp]
	
	\begin{center}
    \begin{tikzpicture}[scale=1,cap=round,>=latex]
    
     \path[draw,-,thick] (0,0) -- (9.5,0);
   
	\foreach \x in {0,1,...,19 } {
       \fill[draw=blue,fill=blue!20!,thick] (\x/2 ,0) circle(1.5pt);
    }

	\foreach \x in {1,5,8,10,13,14,15,16,17 } {
       \fill[draw=blue,fill=red,thick] (\x/2 ,0) circle(2.5pt);
    }

    \end{tikzpicture}
	\end{center}
	\caption{A forbidden configuration.}
        \label{fig-2}
	\end{figure}

Throughout this paper, a necklace is understood to be cyclically symmetric
(unlabelled) formed by diamonds of two colors (i.e., binary).

A {\em configuration} or {\em broken necklace} is one resulting from one of the 
$r+1$ cuts to the left, right or in between the medallion (as described in the 
motivational introduction). 
Figure \ref{fig-1} shows a configuration and Figure \ref{fig-2} 
depicts a forbidden cut.

Label $n$ vertices as $\{1, \, 2, \cdots, n-1, \, n \}$. The {\em 
neighbors} of the vertex $i$ are $i-1$ and $i+1$ for $2 \leq i \leq n-1$; the 
single vertex $2$ for $i=1$ and the single vertex $n-1$ for $i=n$. 
Configurations consist of a linear array of $n$ vertices, $k$ of which are 
{\em marked} or {\em painted red}.  The marked vertices are 
either {\em isolated}, that is, its 
neighbors are not marked or {\em connected}, that is, the sequence of 
vertices $\{ i, i+1, i+2, \cdots, j\}$ are all marked. In the latter case, it 
must be the case that $i=1$ or $j=n$; that is, connected marked vertices 
contain $1$ or $n$. 

\smallskip

\noindent
{\bf Question 1}. Determine the number $Z_{k}(n)$ of configurations up 
to symmetry. 

\smallskip

The problem above, sans restriction, may be interpreted as
a periodic chain made of two kinds of beads. 
The classical result on counting all necklaces with $n$ beads is given 
by MacMahon formula
\begin{equation}
N(n) = \frac{1}{n} \sum_{d | n} \varphi(d) 2^{n/d},
\label{mac}
\end{equation}
\noindent
where the summation runs through all divisors $d$ of $n$, and $\varphi(d)$ is 
the {\em Euler totient} function counting the numbers $1, \, 2, \ldots, \, d$ 
relatively prime to $d$. 


\smallskip

Onofri et al \cite{ovw} introduced a supersymmetric quantum mechanical 
model for a 
system whose degrees of freedom are bosonic and fermionic (creation and 
destruction) operator matrices. In the Hilbert space of this model, the vectors 
can be put in one-to-one correspondence with binary 
necklaces (or periodic linear sequences), where $1$ and $0$ represent 
fermionic and bosonic matrices respectively. {\emph{Pauli's exclusion 
principle}}  provides a Fermi statistics which projects out a subset among 
all necklaces. As such the resulting space allows a genuine depiction 
of supersymmetry. 

The Hilbert space of the above model is generated by states due to single 
trace operators. Since trace is cyclically symmetric, all length $n$ 
necklaces related by cyclic shifts are identified as having the same state 
of $n$ quanta. The Pauli principle is reflected in the fact that the 
fermionic operators are Grassmannian variables, hence they 
anti-commute;  while their bosonic counterparts commute freely as 
scalars. 
The number of necklaces is enumerated by MacMahon's 
formula \eqref{mac}. Some 
of these are excluded as a result of anti-symmetry of planar states, 
dictated by supersymmetry. The following serves as an illustrative example 
with $n=4$. There is a total of $6$ necklaces listed here as linear 
sequences of period $4$, namely
$$0000, \qquad 0001, \qquad 0011, \qquad 0101, \qquad 0111, \qquad 1111.$$
Supersymmetry breaks up this family into forbidden and allowed necklaces as 
follows. Given a sequence, start shifting a digit from right to 
left (prefix to suffix) and repeat until the original sequence is 
recovered. Every time a fermionic operator (that is, a $1$) crosses another 
then a sign change must be registered because of  anti-commuting; each 
bosonic operator (i.e., a $0$) shifts around without any effect. At the end 
of the procedure, if the sequence becomes its own negative then 
call it \it forbidden, \rm otherwise it is \it allowed.  \rm One of the 
main results of \cite{ovw} predicts
\begin{equation}
N_{allowed}(n) = \frac{1}{n} \sum_{\stackrel{d | n}{d \text{ odd}}} 
 \varphi(d) 2^{n/d},
\text{ and } 
N_{forbidden}(n) = \frac{1}{n} \sum_{\stackrel{d | n}{d \text{ even}}} 
 \varphi(d) 2^{n/d}.
\end{equation}

 Let us take a look at each of the above sequences now: $0000$ stays the 
same; $0001\rightarrow 1000\rightarrow 0001$; $0011\rightarrow -1001\rightarrow +1100\rightarrow 0011; 0101\rightarrow -1010\rightarrow -0101; 0111\rightarrow +1011\rightarrow +1101\rightarrow +1110\rightarrow 0111; 1111\rightarrow -1111$. 
Therefore there are four allowed 
$\{0000, 0001, 0011, 0111\}$ and two forbidden $\{0101, 1111\}$ sets of 
necklaces. This agrees with
\begin{equation}
N_{\text{allowed}}=\frac14[\varphi(1)2^4]=4, \quad \text{and} \quad
N_{\text{forbidden}}=\frac14[\varphi(2)2^2+\varphi(4)2^1]=2.
\nonumber
\end{equation}
The interested reader should find the complete story in \cite{ovw}.

\medskip

The central object of the work presented here is a sequence of numbers labeled 
\emph{necklace binomial coefficients} 
$\binom{t}{k}_{\mathfrak{N}}$. These coefficients have properties similar to 
the usual binomial coefficients: recurrences (Corollary \ref{cor-pascal}), 
symmetries  (Corollary \ref{thm-symm}) and an 
explicit formula in terms of the binomial
coefficients (Theorem \ref{form-beta}). Section \ref{sec-polyno} 
discusses a surprising 
result on the zeros of the
necklace polynomial $N_{t}(y)$ (this is the generating function of 
$\binom{t}{k}_{\mathfrak{N}}$). It turns out that all its zeros are inside an
elliptic curve. Moreover, this curve is the same for all values of $t$. 
Section \ref{sec-circ} presents arithmetical and geometrical properties of 
these polynomials via generating function methods.
The necklace binomial coefficients also appear in a physical 
model of Onofri et al \cite{ovw} and in the 
study of symmetries of paraffin molecules in \cite{losanitsch}.

\section{The number of configurations} \label{sec-number}
\setcounter{equation}{0}

In this section the counting problem from the Introduction is rephrased 
and solved. The current format as well as the original formulation 
will be used interchangeably: \\

\noindent
{\em determine the number $Z_{k}(n)$ of 
painting $k$ points in red from a linear array 
of $n$ of them, with the condition that consecutive red points can only 
appear at the beginning or at end of the array. Moreover, arrays that 
are reflections of each other 
should be counted only once. }

\medskip

In order to determine the number of configurations $Z_{k}(n)$ it is 
convenient to begin with a simpler count. 

\begin{Prop}
Let $f_{k}(n)$ be the number of arrangements of $n$ 
vertices with $k$ marked vertices, no consecutive marked ones where
reflections are not 
identified. Then 
\begin{equation}
f_{k}(n) = \binom{n-k+1}{k}.
\end{equation}
\end{Prop}
\begin{proof}
Each such arrangement can be obtained by placing the $k$ marked vertices and 
choosing $k-1$ places to separate them. The count is obtained by eliminating 
the separating spaces.
\end{proof}

\noindent
{\em Reduced  configurations}. The next step is to count those configurations
obtained by cutting the necklace exactly on one side of the medallion. These 
produce linear arrays where clustered vertices appear either at the 
beginning or at the end of the array. Invoking symmetry, only those 
with the medallion at the 
left will be considered. 

\begin{definition}
The function $\beta_{k}(n)$ denotes the number of linear arrays obtained 
by cutting a necklace with $n$ vertices and $k$ marked vertices with a 
medallion at the left of the array. 
\end{definition}

Theorem \ref{form-beta} provides an expression for the function 
$\beta_{k}(n)$ and Theorem \ref{form-alpha} provides a formula for 
$Z_{k}(n)$.

\smallskip

\begin{definition}
Let $g_{k}(n)$ be the number of arrangements of $n$ vertices
with $k$ marked points, no two being consecutively marked and 
identifying symmetric  pairs. 
\end{definition}

\begin{Example}
A numerical example of $g_{k}(n)$ is given 
here. Take  $n=4$ and $k=2$. From the pairs 
$\{ 12, \, 13, \, 14, \, 23, \, 24, \, 34 \}$ eliminate 
$\{ 12, \, 23, \, 34 \}$ for being consecutive somewhere. This leaves 
$\{ 13, \, 14, \, 24 \}$. The pairs are now considered modulo $5$, so that 
$24$ is identified with $13$ (the same as $31$). The final 
allowed list is $\{ 13, \, 14 \}$ 
showing that $g_{2}(4) = 2$. 
\end{Example}

\begin{Thm}
The function $\beta_{k}(n)$ satisfies
\begin{equation}
\beta_{k}(n) = g_{k}(n) + \sum_{r=2}^{k} f_{k-r}(n-r-1).
\end{equation}
\end{Thm}
\begin{proof}
Separate the different configurations into two groups: those with no 
consecutive marked points and those with at least two consecutive ones that 
are marked. The first type is counted by $g_{k}(n)$. Observe that if a 
certain arrangement has two or more adjacent marked vertices, then the 
remaining marked ones have no restrictions due to symmetry. In other words, 
reflection only imposes limitations if the configuration has no adjacent 
marked vertices in it. 

The number of possible consecutive marked points is given by the size of the 
medallion. If this size is $r$, with $2 \leq r \leq k$, then drop $r+1$ places 
from the configuration ($r$ for the medallion and one more at the right-end of it). 
This leaves a total of $n-r-1$ spaces where to place $k-r$ marked vertices.
\end{proof}

The next step is the enumeration of $g_{k}(n)$. This group is divided into 
three disjoint subclasses, those with 
(1) both ends are marked, (2) both ends are unmarked 
and (3) only the left end is marked. In the first class drop the vertices 
at positions $1, \, 2, \, n-1$ and $n$ and observe that the remaining $n-4$
vertices have $k-2$ marked ones and no further restrictions. Therefore there 
are $g_{k-2}(n-4)$ such arrangements. Similarly, the class (2) has 
$g_{k}(n-2)$ elements. Finally, in class (3), drop the first two vertices 
and the last one that is not marked. The remaining $n-3$ vertices have no 
symmetry restriction. The latter are counted by $f_{k-1}(n-3) = 
\binom{n-k-1}{k-1}$ such arrangements. This gives the relation
\begin{equation}
g_{k}(n) = g_{k}(n-2) + g_{k-2}(n-4) + \binom{n-k-1}{k-1}.
\end{equation}

\begin{Thm}
\label{thm-recu-neck}
Let $n = m + 2k-1$ and define $\bar{g}_{k}(m) := g_{k}(m+2k-1)$. Then 
$\bar{g}_{k}$ satisfies 
\begin{equation}
\bar{g}_{k}(m) = \bar{g}_{k-2}(m)  + \bar{g}_{k}(m-2) + \binom{m+k-2}{k-1}.
\label{rec-1}
\end{equation}
\end{Thm}
\begin{proof}
Observe that any valid arrangement counted by $g_{k}(n)$ must satisfy
$n \geq 2k-1$. The rest is elementary.
\end{proof}

The next result was obtained from experimental data generated by
(\ref{rec-1}).

\begin{Example}
The function $\bar{g}_{k}(m)$ is computed for $0 \leq m \leq 3$:
\begin{equation}
\bar{g}_{k}(0) = 1, \, 
\bar{g}_{k}(1) = \left\lfloor \frac{k+2}{2} \right\rfloor, \, 
\bar{g}_{k}(2) = \left\lfloor \frac{(k+2)^{2}}{4} \right\rfloor
\end{equation}
\noindent
and 
\begin{equation}
\bar{g}_{k}(3) =   \sum_{j=0}^{k} (-1)^{k-j} 
\left\{ \sum_{i=0}^{j} \left\lfloor \frac{j+2}{2} \right\rfloor + 
\binom{j+1}{2} \right\}.
\end{equation}
\end{Example}

\begin{Definition}
The relation (\ref{rec-1}) attains a cleaner form by introducing the 
{\em necklace binomial coefficients}
\begin{equation}
\binom{t}{k}_{\mathfrak{N}} := 
\begin{cases}
g_{k}(t+k-1) & \quad \text{ for } 0 \leq k \leq t \\
0 & \quad \text{ otherwise}.
\end{cases}
\end{equation}
\end{Definition}

The next result is a restatement of Theorem  \ref{thm-recu-neck}.

\begin{Cor}
\label{cor-pascal}
The necklace binomial coefficient satisfies the Pascal-type relation
\begin{equation}
\binom{t}{k}_{\mathfrak{N}}  = 
\binom{t-2}{k-2}_{\mathfrak{N}} + 
\binom{t-2}{k-1}+ 
\binom{t-2}{k}_{\mathfrak{N}}.
\label{nb-rec1}
\end{equation}
\end{Cor}

The evaluation of the necklace binomial coefficients is now easy to guess 
and establish using (\ref{nb-rec1}). 

\begin{Thm}
\label{form-beta}
For $0 \leq k \leq t$, it holds that
\begin{equation}
\label{formula-bin}
\binom{t}{k}_{\mathfrak{N}} =  \frac{1}{2}
\begin{cases}
\binom{t}{k} & \quad \text{ for }t \text{ even and }k \text{ odd}, \\
\binom{t}{k} + \binom{\lfloor{t/2 \rfloor}}{\lfloor{k/2 \rfloor} } & 
\quad \text{ otherwise.}
\end{cases}
\end{equation}
\noindent
Moreover, 
\begin{equation}
\beta_{k}(t) = \binom{t-k+1}{k}_{\mathfrak{N}} + 
\sum_{r=2}^{k} \binom{t-k}{r-2}.
\end{equation}
\end{Thm}

Table \ref{table-neck} shows the values of the necklace coefficients:

\medskip

\begin{center}
\begin{tabular}{|c||c|c|c|c|c|c|c|c|c|c|c||}
\hline
$t/k$ & 0 & 1 & 2 & 3 & 4 & 5 & 6 & 7 & 8 & 9 & 10\\
\hline \hline 
1 & 1 & 1 &  &  &  &   &&& &&\\
2 & 1 & 1 & 1 &   &   &  &&&  && \\
3 & 1 & 2 & 2 & 1 &   &  &&&  && \\
4 & 1 & 2 & 4 & 2 & 1  &   &&& && \\
5 & 1 & 3 & 6 & 6 & 3  & 1  &&& && \\
6 & 1 & 3 & 9 & 10 & 9  & 3  & 1 && && \\
7 & 1 & 4 & 12 & 19 & 19  & 12  & 4  & 1 & && \\
8 & 1 & 4 & 16 & 28 & 38  & 28  & 16  & 4 & 1 && \\
9 & 1 & 5 & 20 & 44 & 66  & 66  & 44  & 20 & 5  & 1 & \\
10 & 1 & 5 & 25 & 60 & 110  & 126  & 110  & 60 & 25  & 5 & 1 \\
\hline
\end{tabular}
\label{table-neck}
\end{center}

This array is known as the {\emph{ Losanitsch's triangle}} and information 
about it can be found in Entry A034851 of Neil Sloane's Encyclopedia 
of Integer Sequences. 

\medskip

A series of elementary consequences of (\ref{formula-bin}) are 
presented next. 

\begin{Cor}
The row-sum identity
\begin{equation}
\sum_{k=0}^{t} \binom{t}{k}_{\mathfrak{N}} 
= 2^{t-1} + 2^{\lfloor{(t-1)/2 \rfloor}}
\end{equation}
\noindent
holds.
\end{Cor}

The next statements employ the {\em Fibonacci numbers} $F_{n}$, defined by the 
relation $F_{n} = F_{n-1} + F_{n-2}$ with initial conditions $F_{0}=F_{1} = 1$
and the {\em Lucas numbers} $L_{n}$ defined by the same recurrence 
and with initial conditions $L_{0}=2, \, L_{1} = 1$. 

\begin{Cor}
Let $F_{n}$ and $L_{n}$ as above. Denote
$\tilde{t} := \lfloor{t/2 \rfloor} + 2 + (-1)^{t+1}$. Then
\begin{equation}
\sum_{k=0}^{t} \beta_{k}(t) = \frac{1}{2}( L_{t+2} + F_{\tilde{t}}) -1.
\end{equation}
\end{Cor}

\begin{Cor}
The generating functions
\begin{eqnarray}
\sum_{k=0}^{t} \binom{t}{k}_{\mathfrak{N}} y^{k}  & = & 
\frac{1}{2}(1+y)^{t} + \frac{1}{2} (1 + y^{2})^{\lfloor{t/2\rfloor}} 
(1+y)^{t \bmod 2},  \label{gen-1} \\
\sum_{t \geq 0} \binom{t}{k}_{\mathfrak{N}}x^{t}  & = & 
\frac{(1+x)^{\lfloor{(k+1)/2 \rfloor}} + (1-x)^{\lfloor{(k+1)/2\rfloor}}}
{2(1-x)^{\lceil{(k+1)/2 \rceil}} (1-x^{2})^{\lfloor{(k+1)/2\rfloor}}},
\label{gen2}
\end{eqnarray}
\noindent
and
\begin{equation}
\sum_{t,k \geq 0} \binom{t}{k}_{\mathfrak{N}} x^{t}y^{k} =
\frac{1}{2(1-x-y)} + \frac{2+x}{2(1-x^{2}-y)},
\label{gen-3}
\end{equation}
\noindent
hold.
\end{Cor}

\begin{Cor}
\label{thm-symm}
The necklace binomial coefficients are symmetric, that is,
\begin{equation}
\binom{t}{k}_{\mathfrak{N}}  = \binom{t}{t-k}_{\mathfrak{N}}
\label{neck-symm}
\end{equation}
\noindent
for $0 \leq k \leq t$. 
\end{Cor}

\begin{Cor}
The function $\bar{g}$ is symmetric; that is,
\begin{equation}
\bar{g}_{k}(m) = \bar{g}_{m}(k).
\label{sym-bin}
\end{equation}
\end{Cor}

\begin{proof}
This is a restatement of (\ref{neck-symm}). An 
alternative proof of the symmetry (\ref{sym-bin}) is obtained from the 
recurrence (\ref{rec-1}). Simply express it in two different forms
\begin{eqnarray}
\bar{g}_{k}(m) - \bar{g}_{k-2}(m) & = & \bar{g}_{k}(m-2) + \binom{m+k-2}{k-1}, 
\label{rec1-a} \\
\bar{g}_{k}(m) - \bar{g}_{k}(m-2) & = & \bar{g}_{k-2}(m) + \binom{m+k-2}{k-1}.
\nonumber
\end{eqnarray}
\noindent
The result now follows by induction and the symmetry of the 
binomial coefficients.
\end{proof}

The next theorem provides a combinatorial proof of the symmetry rule 
(\ref{sym-bin}).

\begin{Thm}
The symmetry $\bar{g}_k(m)=\bar{g}_m(k)$ holds.
\end{Thm}
\begin{proof}
The assertion amounts to $g_k(m+2k-1)=g_m(k+2m-1)$. Take a linear array of 
$n$ nodes and its $2$-coloring (red $r$ or white $w$). By definition, $g_k(n)$ 
enumerates all possible ways of coloring $k$ nodes in red with the rule: 
(1) no two reds are consecutive; (2) two such arrays are equivalent if they 
relate by reflection. According to (1), it must be that the first $k-1$ 
reds are each followed by white. Thus, any selection of $k$ reds can be 
interpreted as choosing the $(k-1)$ pairs $rw$ and a free $r$. For each pair 
$rw$, trim-off the $w$ as well as its sitting node. That means, when $n=m+2k-1$ 
then the number of nodes reduces to $m+k$ and hence $g_k(m+2k-1)$ induces an 
equivalent counting of $(m+k)$-nodes of which $k$ are red (note: rule (1) is 
absent but rule (2) stays). Similarly, $g_m(k+2m-1)$ tantamount to the 
counting of $(m+k)$-nodes of which $m$ are white. But, it is obvious that 
coloring $k$ nodes red on an $(m+k)$-array is equivalent to the coloring 
of $m$ nodes in white. This gives the required bijection. The proof is 
complete.
\end{proof}

\begin{Example}
This example demonstrates the above proof; i.e. $g_k(m+2k-1)=g_m(k+2m-1)$. 
Take $m=2$ and $k=3$. Then, $g_3(7)$ and $g_2(6)$ count respectively the 
cardinality of sets
\begin{equation}
A:=\{rwrwrww, rwrwwrw, rwrwwwr, rwwrwrw, rwwrwwr, wrwrwrw \} 
\nonumber 
\end{equation}
\noindent
and 
\begin{equation}
B:\{rwrwww, rwwrww, wrwrww, rwwwrw, wrwwrw, rwwwwr\}.
\nonumber
\end{equation}
The set $B$ after color-swapping turns to
$$B_1:=\{wrwrrr, wrrwrr, rwrwrr, wrrrwr, rwrrwr, wrrrrw\}.$$
The two sets $A$ and $B_1$ are now mapped ($w$-trimmed and 
$r$-trimmed, respectively) to
\begin{equation}
A_1:=\{rrrww, rrwrw, rrwwr, rwrrw, rwrwr, wrrrw\},
\nonumber
\end{equation}
\noindent
and 
\begin{equation}
B_{11}:=\{wwrrr, wrwrr, rrwrr, wrrwr, rwrwr, wrrrw\}.
\nonumber
\end{equation}

The bijection between $A_1$ and $B_{11}$ is clearly exposed; that is, reflect 
$B_{11}$ to get the set
\begin{equation}
B_{111}:=\{rrrww, rrwrw, rrwrr, rwrrw, rwrwr, wrrrw\}.
\nonumber
\end{equation}
\end{Example}

The full counting solution to the configuration problem is presented next. 

\begin{Thm}
\label{form-alpha}
The total number $Z_{k}(t)$ of possible linear configurations of $k$ 
diamonds (with or without a medallion) on $t$ nodes is given by 

\begin{equation}
Z_{k}(t) =  
\sum_{j \geq 0} \binom{t-k-1}{k-2j}_{\mathfrak{N}} 
+ \sum_{j \geq 0} \left\lfloor \frac{j+1}{2} \right\rfloor \binom{t-k-1}{k-j}.
\nonumber
\end{equation}
\end{Thm}
\begin{proof}
Catalog the diamonds according to whether they are: (1) an equal number of 
clusters; (2) unequal number of clustered diamonds on the two end-nodes. 
However many are remaining to be mounted in the interior, case (1) is 
affected by the reflection but those in case (2) are not. It follows that 
the first case is enumerated by the function $g_{k}(t)$ (equivalently, by 
necklace binomials) while the function $f_{k}(t)$ is the right choice 
for the second category. The details are omitted. 
\end{proof}

The necklace coefficients are given as 
Entry A005994 in Neil Sloane's Encyclopedia
of Integer Sequences. The reader will find there information on the 
connection between
 $\binom{t}{k}_{\mathfrak{N}}$ and the so-called {\em paraffin numbers}. 
The chemist S. M. Losanitsch studied in \cite{losanitsch} the so-called 
{\em alkane numbers} (called here the necklace numbers) in his investigation
of symmetries manifested by rows of paraffin (hydrocarbons). In the molecule
of an {\em alkane} (also known as a paraffin), for $n$ carbon atoms there 
are $2n+2$ hydrogen atoms (i.e. the form $C_{n}H_{2n+2}$). Each carbon atom $C$ 
is linked to four other atoms (either $C$ of $H$); each hydrogen atom is 
joined to one carbon atom. The figures in the Appendix show all possible alkane 
bonds for $1 \leq n \leq 5$. There are $1, \, 1, \, 1, \, 2, \, 3$ possible 
alignments, respectively.

\medskip

\noindent
{\bf A geometric interpretation}. Given a finite group $G$, it is a classical 
problem to find the generators of the ring of polynomial invariants under the 
action of $G$. The {\em Molien series} $M(z;G)$ is the generating function 
that counts the number of linearly independent homogeneous polynomials of a given total degree $d$ that are 
invariants for $G$. It is given by 
\begin{equation}
M(z;G) = \frac{1}{|G|} \sum_{g \in G} \frac{1}{\text{det}(I - zg)} =
\sum_{i=0}^{\infty} b_{i}z^{i}.
\end{equation}
\noindent
Thus, the coefficients $b_{i}$ record the number of linearly independent 
polynomials of total degree $i$. 

Now assume $k = 2m-1$. Then (\ref{gen2}) becomes
\begin{equation}
\sum_{i \geq 0} \binom{i+2m-1}{2m-1}_{\mathfrak{N}} z^{i} = 
\frac{1}{2} \frac{1}{(1-z)^{2m}} + \frac{1}{2} \frac{1}{(1-z^{2})^{m}}.
\end{equation}
\noindent
This is recognized as 
\begin{equation}
\frac{1}{2} \frac{1}{(1-z)^{2m}} + \frac{1}{2} \frac{1}{(1-z^{2})^{m}} = 
\frac{1}{|G|} \sum_{g \in G} \frac{1}{\text{det}({\mathbf{1}}_{2m}-zg)},
\end{equation}
\noindent
where $G$ is the symmetric group $S_{2}$ and the summation runs through the 
$2m$-dimensional group representation of the elements $g$ in 
$GL_{2m}(\mathbb{C})$. The argument below shows that the series is indeed a
Molien series for the ring of invariants under the action of $S_{2}$. More 
specifically, the ring of invariants under consideration is 
$\mathbb{C}[X;Y]^{{\mathbb{s}}_{2}}$ where $X = (x_{1},\ldots,x_{m})$ 
and $Y = (y_{1},\ldots,y_{m})$. The action is given by $x_{l} \mapsto y_{l}$ 
for $l=1, \ldots, n$. 

Let $\sigma$ be the matrix 
$\sigma = \begin{pmatrix} 0 & 1 \\ 1 & 0 \end{pmatrix}$ and let $\pi$ 
be the tensor product $\pi = \sigma \otimes {\mathbf{1}}_{m}$ resulting
in a $2m \times 2m$ matrix which has four 
blocks of size $m \times m$ with the off-diagonal blocks being the identity 
matrix and the diagonals blocks being zero. The matrix group 
generated by $\pi$ in $GL_{2m}$ is $S_{2}$. Consequently,
\begin{equation}
\text{det}( \mathbf{1}_{2m} - z \pi^{2}) = 
\text{det}( \mathbf{1}_{2m} - z {\mathbf{1}}_{2m}) = (1-z)^{2m}
\end{equation}
\noindent
and $\text{det}( \mathbf{1}_{2m} - z \pi) = \text{det}(\rho \otimes 
{\mathbf{1}}_{2m} )$, with 
$\rho = \begin{pmatrix} 1 & -z \\ -z & 1 \end{pmatrix}$. Since 
$\text{det}(A \otimes B) = \text{det}(A)^{m} \, \text{det}(B)^{m}$, it must 
be that $\text{det}( {\mathbf{1}}_{2m} - z \pi) = (1-z^{2})^{m}$.  \\

These observations are summarized in the next statement.

\begin{Thm}
Consider the action of $\mathbb{Z}_{2}$ on $\mathbb{C}[x_{1}, \cdots,
x_{m},y_{1}, \cdots, y_{m} ]$ given by $x_{l} \mapsto y_{l}$. Then, the number 
of linearly independent invariant polynomials of total degree $i$ is given
by the necklace binomial coefficient
$\binom{i+2m-1}{2m-1}_{\mathfrak{N}}$.
\end{Thm}

\section{The necklace polynomials} \label{sec-polyno}
\setcounter{equation}{0}

In this section we discuss properties of the {\em necklace polynomials}
 defined by 
\begin{equation}
N_{t}(y) = \sum_{k=0}^{t} \binom{t}{k}_{\mathfrak{N}}y^{k}.
\end{equation}

The explicit formula 
\begin{equation}
N_{t}(y) = \frac{1}{2} (1+y)^{t} + \frac{1}{2}(1+y^{2})^{\lfloor{t/2 \rfloor}}
(1+ y)^{t \, \bmod 2}
\label{Thm-explicit}
\end{equation}
\noindent
is given in \eqref{gen-1}. 

\begin{Example}
The first few values of $N_{t}(y)$ are given by 
\begin{eqnarray}
N_{1}(y) & = & 1 + y \nonumber  \\
N_{2}(y) & = & 1 + y + y^{2} \nonumber \\
N_{3}(y) & = & N_{1}(y)N_{2}(y) \nonumber \\
N_{4}(y) & = & 1 + 2y + 4y^{2} + 2y^{3} + y^{4} \nonumber \\
N_{5}(y) & = & N_{1}(y)N_{4}(y) \nonumber \\
N_{6}(y) & = & N_{2}(y) (1 + 2y + 6y^{2} + 2y^{3} + y^{4}) \nonumber \\
N_{7}(y) & = & N_{1}(y) N_{2}(y) (1 + 2y + 6y^{2} + 2y^{3} + y^{4}) \nonumber \\
N_{8}(y) & = & 1 + 4y + 16y^{2} + 28y^{3} + 38y^{4} + 28y^{5} + 
16y^{6} + 4y^{7} + y^{8}. \nonumber
\end{eqnarray}
\end{Example}

The sequence of necklace polynomials have some interesting divisibility 
properties. The results presented below began with the empirical observation
that, for $t$ odd, $N_{t}(y) = N_{1}(t) N_{t-1}(y)$. 

\begin{Cor}
Let $j \in \mathbb{N}$ and $t \in \mathbb{N}$. Then $N_{j}(y)$ divides 
$N_{(2t-1)j}(y)$. 
\end{Cor}
\begin{proof}
This is a direct consequence of the explicit formula given in Theorem 
\ref{Thm-explicit}.
\end{proof}

\begin{Problem}
Prove that $N_{2^{j}}(y)$ is irreducible.
\end{Problem}

\medskip

Many polynomials appearing in combinatorics are {\em unimodal}; that is, 
there is an index $n^{*}$ such that the coefficients increase up to 
$n^{*}$ and decrease from that point on. A stronger property is that of 
{\em logconcavity}: the polynomial $P(x) = \sum_{k=0}^{n}a_{k}x^{k}$ is 
logconcave if $a_{k}^{2} - a_{k-1}a_{k+1} \geq 0 $ for $1 \leq k \leq n-1$. 
The reader is referred to 
\cite{brenti, stanley} for surveys on these issues. 

\smallskip

The explicit expression (\ref{formula-bin}) gives an 
elementary proof of the next statement. 

\begin{Thm}
The necklace binomial coefficients are unimodal.
\end{Thm}
\begin{proof}
The inequality 
\begin{equation}
\binom{t}{k}_{\mathfrak{N}} \leq \binom{t}{k+1}_{\mathfrak{N}} 
\end{equation}
\noindent
for $0 \leq k \leq \lfloor{ t/2 \rfloor}$ and the symmetry of the 
necklace binomial coefficients, established in Theorem \ref{thm-symm},
give the result.
\end{proof}

\begin{Thm}
The polynomial $N_{t}(y)$ is logconcave. 
\end{Thm}
\begin{proof}
Use (\ref{formula-bin}) and separate cases according to the parity of 
$t$ and $k$.
\end{proof}

\begin{Problem}
Let $\mathfrak{L} \{ a_{n} \} := \{ a_{n}^{2} - a_{n-1}a_{n+1} \}$ be an 
operator defined on nonnegative sequences. Therefore, a 
polynomial $P(x)$ is logconcave if 
$\mathfrak{L}$ maps its coefficients into a 
nonnegative sequence. The polynomial $P$ is 
called $k$-logconcave if $\mathfrak{L}^{(j)}(P)$ is nonnegative for 
$0 \leq j \leq 
k$. A sequence is called {\em infinitely logconcave} if it is $k$-logconcave 
for every $k \in \mathbb{N}$. 

A recent result of P. Br\"{a}nd\'{e}n \cite{branden} 
proves that if a polynomial 
$P$ has only real and negative zeros, then the sequence of its coefficients 
is infinitely logconcave. The sequence of 
binomial coefficients satisfies this property.

\smallskip

The question proposed here is to prove that $N_{t}(y)$ is 
infinitely logconcave.
\end{Problem}

\medskip

There is a well-established connection between unimodality questions and 
the location of the zeros of a polynomial. For example, a polynomial with 
all its zeros real and negative is logconcave \cite{wilf}. This motivated 
the computation of the zeros of $N_{t}(y)$. Figure \ref{fig-3} shows the 
zeros of $N_{100}(y)$.


{{
\begin{figure}[ht]
\begin{center}
\centerline{\epsfig{file=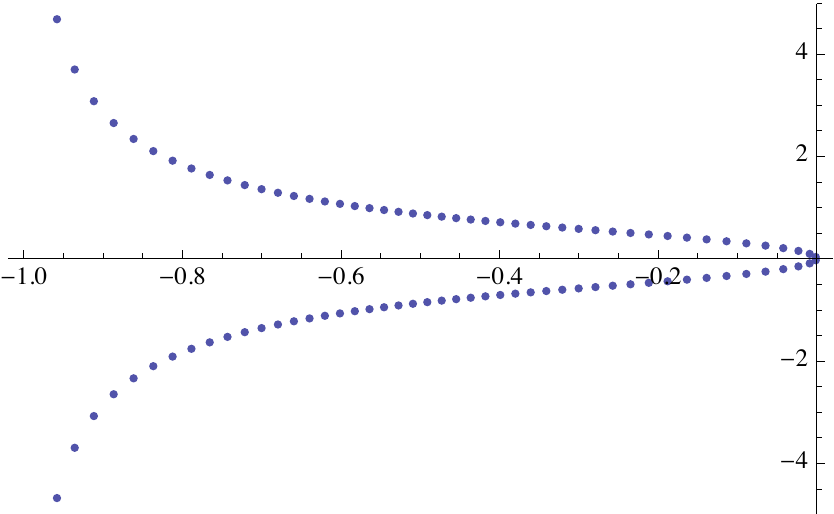,width=20em,angle=0}}

\caption{The zeros of the necklace polynomial $N_{100}(y)$.}
\label{fig-3}
\end{center}
\end{figure}
}}



\medskip

\begin{Thm}
Let $y = a + ib$ be a root of the necklace polynomial $N_{t}(y) = 0$. For 
$a  \neq -1$, define 
the new coordinates $u = 1/(1+a)$ and $v = b/(1+a)$. Then $(u,v)$ is on the 
elliptic curve $v^{2} = u^{3}-2u^{2}+2u-1$. 
\end{Thm}
\begin{proof}
Any zero of $N_{t}(y)$ satisfies
\begin{equation}
(1+y)^{t} = - \begin{cases}
              (1+y^{2})^{t/2} & \quad \text{ if } t \text{ is even} \\
              (1+y^{2})^{(t-1)/2}(1+y) & \quad \text{ if } t \text{ is odd.}
              \end{cases}
\end{equation}
\noindent
Taking the complex modulus 
produces $|1 + y|^{4} = | 1 + y^{2}|^{2}$. In terms of $y = a + ib$ this 
equation becomes
\begin{equation}
b^{2} = - \frac{a(a^{2}+a+1)}{1+a}.
\label{curve-1}
\end{equation}
\noindent
The transformation $1+a  = 1/u$ and $b = v/u$ leads to equation
\begin{equation}
v^{2} = u^{3} - 2u^{2} + 2u - 1 = (u-1)(u^{2}-u+1),
\label{cubic}
\end{equation}
\noindent
as claimed.
\end{proof}

\begin{Note}
The collection of points on an elliptic curve $\mathfrak{E}$, such as 
(\ref{cubic}), has been the subject of research since the $18$th century. 
The general equation of such a curve is written as
\begin{equation}
y^{2} + a_{1}y = x^{3} + a_{2}x^{2} + a_{4}x + a_{6}
\end{equation}
\noindent
and if $x,y \in P(\mathbb{C}^{2})$, the complex projective space,  then 
$\mathfrak{E}$ is a torus. The addition
of this torus is expressed on the cubic in a geometric form: to add $P_{1}$
and $P_{2}$, form the line joining them and define $P_{3}:= P_{1} \oplus P_{2}$ 
as the reflection of the third point of intersection of this line with the 
cubic curve. This addition rule is expressed in coordinate form: the 
general formula given in \cite{silverman-I}. Let 
$P_{1} = (x_{1},y_{1})$ and $P_{2} = (x_{2},y_{2})$. Define 
\begin{eqnarray}
\lambda = \begin{cases}
      \frac{y_{2}-y_{1}}{x_{2}-x_{1}} & \quad \text{ if } x_{2} \neq x_{1} 
\nonumber \\
\frac{3x_{1}^{2} - 4x_{1} +2}{2y_{1}} & \quad \text{ if } x_{2} = x_{1},
\end{cases}
& \text{ and } & 
\nu = \begin{cases}
 \frac{y_{1}x_{2}-y_{2}x_{1}}{x_{2}-x_{1}} & \quad \text{ if } x_{2} \neq x_{1} 
\nonumber \\
\frac{-x_{1}^{3} +2x_{1} -2}{2y_{1}} & \quad \text{ if } x_{2} = x_{1}.
\end{cases}
\end{eqnarray}
\noindent
The $P_{3} = (x_{3},y_{3})$ is given by 
\begin{equation}
x_{3}  = \lambda^{2} +2 - x_{1}-x_{2} \text{ and } 
y_{3} = - \lambda x_{3} - \nu. 
\nonumber
\end{equation}

\medskip

Aside from the point $P_{0} = (1,0)$, the table below shows a collection 
of points on the curve $\mathfrak{E}$ 
obtained using Mathematica. The notation
\begin{multline}
\gamma= \sqrt{3 + 2 \sqrt{3}}, \, \delta = \sqrt{\sqrt{5}-2}, \, 
\tau= \sqrt{24 + 14 \sqrt{3}}, \, 
\sigma= 2\sqrt{2(11+ 5 \sqrt{5})}, 
\\
\omega_{1} = 2 + \sqrt{3}, \, \omega_{2} = 2(3 + \sqrt{5}), \, 
\omega_{3} = 3 + 2 \sqrt{3}
\nonumber
\end{multline}
\noindent 
is employed.

\medskip

\begin{table}[ht]
\begin{center}
\begin{tabular}{||c|c|c|c||}
\hline 
\text{Name} & $u$ & $v$ & \text{Root of} $N_{t}(y) = 0$  \\
\hline
$P_{1}$ & $2$ & $- \sqrt{3}$ & 2 \nonumber \\
$P_{2}$ & $2$ & $+ \sqrt{3}$ & 2 \nonumber \\
$P_{3}$ & $\omega_{1} - \gamma$ & 
$\omega_{3}- \tau$ & 6 \nonumber \\
$P_{4}$ & $\omega_{1}- \gamma$ & 
$- \omega_{3} + \tau $ & 6 \nonumber \\
$P_{5}$ & $\omega_{1} + \gamma$ & 
$- \omega_{3} - \tau $ & 6 \nonumber  \\
$P_{6}$ & $\omega_{1} + \gamma$ & 
$ \omega_{3} + \tau $ & 6 \nonumber \\
$P_{7}$ & $(1 + \delta) \omega_{2}$ 
& $ \omega_{2} + \sigma $  & 4 \nonumber \\
$P_{8}$ & $(1 + \delta ) \omega_{2} $ 
& $ - ( \omega_{2}+ \sigma )$  & 4 \nonumber \\
$P_{9}$ & $(1 - \delta ) \omega_{2}$ 
& $ \omega_{2}- \sigma $  & 4 \nonumber \\
$P_{10}$ &$ (1 - \delta ) \omega_{2} $ 
& $ -( \omega_{2}  - \sigma)$  & 4 \nonumber  \\
\hline
\end{tabular}
\end{center}
\caption{Some points on the elliptic curve $\mathfrak{E}$.}
\end{table}

\medskip

The notation {\em necklace point} refers to a point $(u,v)$ on the 
elliptic curve $\mathfrak{E}$ that is produced by the zero $y = a + ib$ of a 
necklace polynomial via 
the transformation $1+a  = 1/u$ and $b = v/u$. 
The addition of two necklace points sometimes yields another one. For instance,
$P_{1} \oplus P_{1} = P_{0}$ and $2P_{3}:= P_{3} \oplus P_{3} = P_{2}$. On the 
other hand, the set of necklace points is not closed under addition:
\begin{equation}
\nonumber
P_{1} \oplus P_{7} = 
\frac{1}{2} \left( 7 + 3 \sqrt{5} + \sqrt{66 + 30 \sqrt{5}}  \, \right) - 
\frac{I}{2} \left( 21 + 9 \sqrt{5} + \sqrt{30(29 + 13 \sqrt{5})} \, \right). 
\end{equation}
\noindent
The minimal polynomial for this number is 
$y^{8} -28 y^{7} + 1948y^{6} - 5236y^{5} + 4858y^{4} - 3988y^{3} 
+ 7156y^{2} - 6040y + 2245$. This polynomial does not divide a $N_{t}(y)$
for $1 \leq t \leq 1000$. It is conjectured that it never does. 
\end{Note}

\begin{Note}
Equation (\ref{curve-1}) shows that any root of $N_{t}(y)$
must satisfy $-1 \leq \realpart{y} \leq 0$. Observe that $y=0$ is never a root. 
\end{Note}

\begin{Note}
The change of variables $u \mapsto u+1$ transforms the curve $\mathfrak{E}$
into the form $v^{2} = u^3 + u^2+u$. This curve appears as $48a4$ in 
Cremona's table of elliptic curves, available at 
\begin{center}
\url{http://www.ma.utexas.edu/users/tornaria/cnt/cremona.html?conductor=48}
\end{center}
\noindent
The discriminant of the cubic is negative. Therefore the curve has a single 
real component. This is seen in Figure \ref{fig-3}.
\end{Note}

\begin{Problem}
The zeros of the polynomial $N_{t}(y)$ are algebraic numbers lying on the 
elliptic curve (\ref{curve-1}). The points on that curve with algebraic 
coordinates form a subgroup $\mathcal{A}$ under the 
addition described above. The question is 
to characterize in $\mathcal{A}$ the set coming from necklace points.
\end{Problem}

\section{Necklaces and their progeny} \label{sec-circ}
\setcounter{equation}{0}

This section explores the enumeration of certain special 
necklaces and their generating functions. The latter is applied to 
the computation of some Molien series. A {\em circuit graph} is a graph 
consisting of $n$ vertices placed on a circle with some of them colored by red.

\begin{Prop}
The total number of
$n$-bead (circular) necklaces on which a red-red string is forbidden 
is given by 
\begin{equation}
W(n) = \frac{1}{n} \sum_{d | n} \varphi \left( \frac{n}{d} \right) L_{d}.
\label{w-for}
\end{equation}
\end{Prop}
\begin{proof}
A standard application of Burnside's lemma.
\end{proof}

\begin{Example}
For $n = p$ prime, formula (\ref{w-for}) gives 
\begin{equation}
W(p) = \frac{(p-1) + L_{p}}{p}.
\end{equation}
\noindent
It follows that $L_{p} \equiv  1 \bmod p$. Similarly, for 
$n = p^{2}$, (\ref{w-for}) gives 
\begin{equation}
p^{2}W(p^{2}) = L_{p^{2}} + (p-1)L_{p} + p(p-1).
\end{equation}
\noindent
It follows that 
\begin{equation}
L_{p^{2}} \equiv L_{p} + 1 \bmod p^{2}.
\end{equation}
\noindent
These are well-known results \cite{hardy4}. 
\end{Example}

A more distinguishing count is provided by defining
$W_{k}(n)$ to be the number of $n$-bead 
necklaces on which a red-red string is forbidden, 
consisting of exactly $k$ red beads.
In order to accomodate the possibility that $k=0$, we define 
$W_{0}(n) :=1$ (this is justifiable since $W_{0}(n) = \frac{1}{n} 
\sum_{d | n} \varphi(d) = 1$) . 

\begin{Thm}
\label{formula-wn}
The function $W_{k}(n)$ is given by 
\begin{equation}
W_{k}(n) = \frac{1}{n-k} \sum_{d | n,k} \varphi(d) 
\binom{ \tfrac{n}{d} - \tfrac{k}{d} }{\tfrac{k}{d}}. 
\end{equation}
\end{Thm}
\begin{proof}
It follows directly from Burnside's lemma.
\end{proof}

\begin{Cor}
The identity 
\begin{equation}
\sum_{k=0}^{\lfloor{n/2 \rfloor}} \frac{1}{n-k} \sum_{d | n,k} \varphi(d) 
\binom{ \tfrac{n}{d} - \tfrac{k}{d} }{\tfrac{k}{d}} = 
\frac{1}{n} \sum_{d | n} \varphi(d) L_{n/d}
\end{equation}
\noindent
holds. 
\end{Cor}
\begin{proof}
The assertion follows from the combinatorial identity
\begin{equation}
\sum_{k \geq 0} W_{k}(n) = W(n).
\end{equation}
\end{proof}

\begin{Thm}
For $n \in \mathbb{N}$ define 
\begin{equation}
V_{d}(x) = \left( \frac{1 - \sqrt{1 + 4x}}{2} \right)^{d} + 
\left( \frac{1 + \sqrt{1 + 4x}}{2} \right)^{d}.
\end{equation}
\noindent
Then the row-sum generating function of $W_{k}(n)$ is given by 
\begin{equation}
F_{n}(x) := \sum_{k=0}^{\lfloor{n/2 \rfloor}} W_{k}(n) x^{k} 
= \frac{1}{n} \sum_{d | n} 
\varphi \left( \frac{n}{d} \right) V_{d} \left( x^{n/d} \right).
\label{row-sum}
\end{equation}
\end{Thm}
\begin{proof}
The proof is based on the identity
\begin{equation}
\frac{1}{m} V_{m}(x) = \sum_{k=0}^{\lfloor{m/2 \rfloor}}  
\frac{1}{m-k} \binom{m-k}{k}x^{k},
\end{equation}
\noindent
which is easy to verify. This is applied to 
\begin{eqnarray}
\sum_{k=0}^{\lfloor{n/2 \rfloor} } W_{k}(n) x^{k} & = & 
\sum_{d | n} \varphi(d) \sum_{k \geq 0} \frac{1}{n - dk} 
\binom{ \tfrac{n}{d} - k}{k} x^{dk} \nonumber \\
& = & 
\sum_{d | n} \frac{\varphi(d)}{d} \, \sum_{k \geq 0} 
\frac{1}{ \tfrac{n}{d} - k } 
\binom{\tfrac{n}{d} - k}{k} x^{dk}.
\nonumber 
\end{eqnarray}
\noindent
The result follows from here. 
\end{proof}

\begin{Example}
For $p$ prime, the polynomial $F_{p}(x)$, defined in \eqref{row-sum}, is given by 
\begin{eqnarray}
F_{p}(x) & = & \sum_{k=0}^{\lfloor{p/2 \rfloor}} 
\frac{1}{p-k} \binom{p-k}{k}x^{k} \nonumber \\
& = & \frac{(p-1)2^{p} (1 - \sqrt{1+4x})^{p} + (1 + \sqrt{1+4x})^{p}}
{p \cdot 2^{p}}.
\nonumber
\end{eqnarray}
\end{Example}

\begin{Example}
Put $n = 3k+1$ in (\ref{formula-wn}) to obtain 
$W_{k}(3k+1) = \frac{1}{2k+1} \binom{2k+1}{k}$, the 
Catalan numbers.
\end{Example}

\begin{Example}
For $n \in \mathbb{N}$, and with $L_{n}$ denoting the Lucas number,
\begin{equation}
\sum_{k=0}^{\lfloor{n/2 \rfloor}} 
\frac{1}{n-k} \binom{n-k}{k} = \frac{1}{n} L_{n}.
\end{equation}
\noindent
This is obtained from setting $x=1$ in (\ref{row-sum}).
\end{Example}

\smallskip

\begin{Thm}
The ordinary generating function  for the 
diagonals of $W_{k}(n)$ is given by 
\begin{equation}
\sum_{n \geq k} W_{k}(n) x^{n} = \frac{1}{k} \sum_{d | k} 
\frac{\varphi(d) x^{2k}}{(1-x^{d})^{k/d}}.
\end{equation}
\noindent
In its lowest terms, the denominator of this rational function takes the form
\begin{equation}
\prod_{d | k} (1-x^{d})^{\varphi(k/d)} = 
\prod_{d | k} \Phi_{d}(x)^{k/d},
\end{equation}
\noindent
where $\Phi_{d}(x)$ is the $d$-th cyclotomic polynomial given in terms of 
the Mobius $\mu$-function as
$\Phi_{d}(x) = \prod_{c|d} (1-x^{d/c})^{\mu(x)}$.
\end{Thm}
\begin{proof}
The result follows from the Taylor series expansion
\begin{equation}
\frac{x^{2m}}{m(1-x)^{m}} = \sum_{j \geq m} \frac{1}{j-m} \binom{j-m}{m} 
x^{j}.
\end{equation}
\end{proof}

\noindent
{\bf A geometric interpretation}. The 
above generating function $\sum_{n \geq k} W_{k}(n)x^{n}$ is the Molien 
series $W(x; {\mathbb{Z}}_{k})$ for the ring of invariants 
$\mathbb{C}[X]^{\mathbb{Z}_{k}}$ where $X = (x_{1},\ldots, x_{k})$. In this 
case, the group $\mathbb{Z}_{k}$ is identified with its $k$-dimensional 
group representation in $GL_{k}(\mathbb{C})$. More concretely, 
$\mathbb{Z}_{k} \cong \langle{ {\mathbf{e}}_{k} \rangle}$ where 
$ {\mathbf{e}}_{k}$ is the $k \times k$ permutation matrix such that 
${\mathbf{e}}[i,j] = 1$ if $j=i+1; {\mathbf{e}}[k,1]=1$ and ${\mathbf{e}}[i,j]
= 0$, otherwise. Let $RP(d)$ be the set of positive integers less than $d$ and 
relatively prime to $d$. Partition the integer interval $[ k ]$
into the disjoint union
\begin{equation}
[k] = \{ 1, \, 2, \, \ldots, \, k \} = \bigcup_{d | k} \frac{k}{d} \, RP(d).
\end{equation}
\noindent
This relation is reminiscent of the well-known identity 
$ k = \sum_{d | k} \varphi(d)$. Then,
\begin{eqnarray}
W(x;{\mathbb{Z}}_{k}) & = & \frac{1}{| {\mathbb{Z}}_{k} |}
\sum_{j=1}^{k} \frac{1}{\text{det}( {\mathbf{1}}_{k} - 
x {\mathbf{e}}_{k}^{j})}  \nonumber \\
& = & \frac{1}{k} \sum_{d | k} \frac{\varphi(d)}{\text{det}( 
\mathbf{1}_{k} - x {\mathbf{e}_{k}^{k/d}} )} \nonumber \\
& = & \frac{1}{k} \sum_{d | k} \frac{\varphi(d)}{\text{det}( 
(\mathbf{1}_{d} - x {\mathbf{e}_{d}}) \otimes {\mathbf{1}}_{k/d})} \nonumber \\
& = & \frac{1}{k} \sum_{d | k} \frac{\varphi(d)}{\text{det}( 
(\mathbf{1}_{d} - x {\mathbf{e}_{d}})^{k/d}} \nonumber \\
& = & \frac{1}{k} \sum_{d | k} \frac{\varphi(d)}{(1-x^{d})^{k/d}}. \nonumber
\end{eqnarray}

\smallskip

These findings are stated in the next result. 

\begin{Prop}
The number of linearly independent homogeneous polynomials, of total degree $n$, for 
the ring of invariants $\mathbb{C}[X]^{\mathbb{Z}_{k}}$ equals 
$$ \frac{1}{n+k} \sum_{d |n,k} \varphi(d) \binom{\tfrac{n}{d} + \tfrac{k}{d}}
{\tfrac{k}{d}}.$$
\end{Prop}

\medskip

\section{A sample of the computation of zeros} \label{sec-comp-zeros}
\setcounter{equation}{0}

Motivated by the interesting properties of the zeros of necklace polynomials,
this section presents some computational graphics showing the zeros of the 
polynomials $F_{n}(x)$. Figure \ref{fig-4}
shows the location of the roots of $F_{1000}(x)$.


{{
\begin{figure}[ht]
\begin{center}
\centerline{\epsfig{file=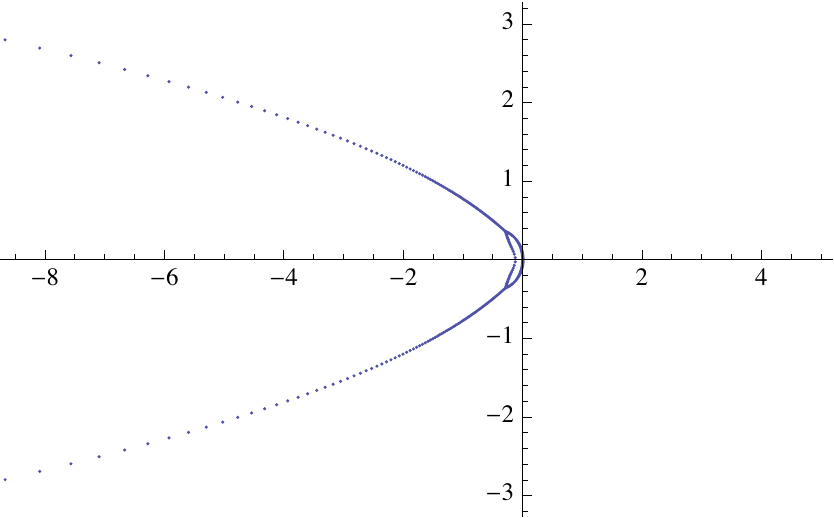,width=30em,angle=0}}

\caption{The zeros of $F_{1000}(x)$.}
\label{fig-4}
\end{center}
\end{figure}
}}


\newpage

The next four 
figures show a selection of regions from the set of the roots 
of all the polynomials $F_{n}(x)$ for $3 \leq n \leq 1000$.  The caption 
indicates the range depicted. 


{{
\begin{figure}[ht]
\begin{minipage}[t]{16em}
\centering
\epsfig{file=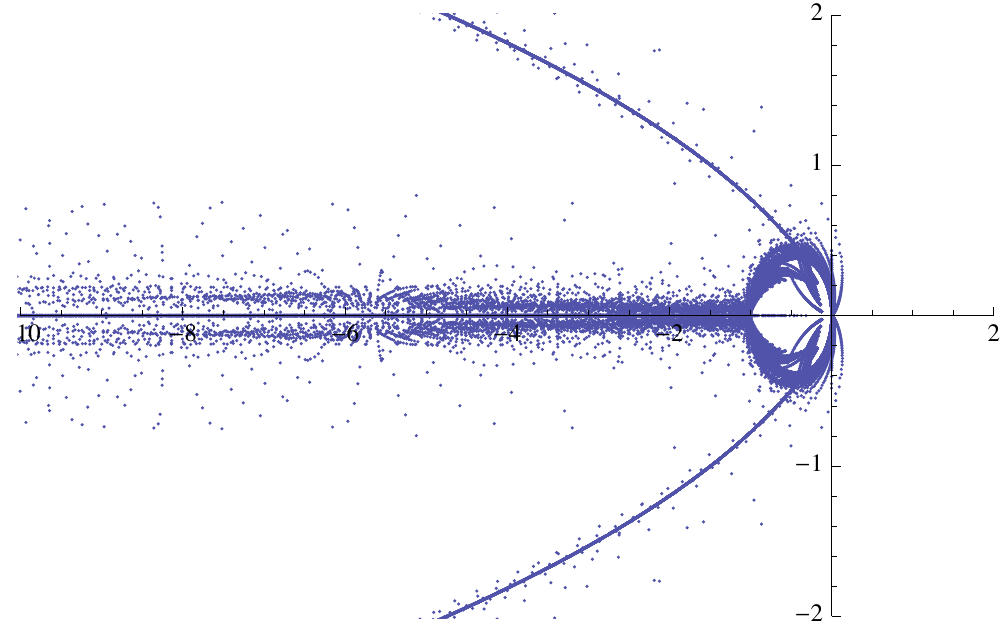,width=15em,angle=0}

\caption{$[-10,2] \times [-2,2]$.}
\label{fig-1-zero}
\end{minipage}%
\begin{minipage}[t]{16em}
\centering
\epsfig{file=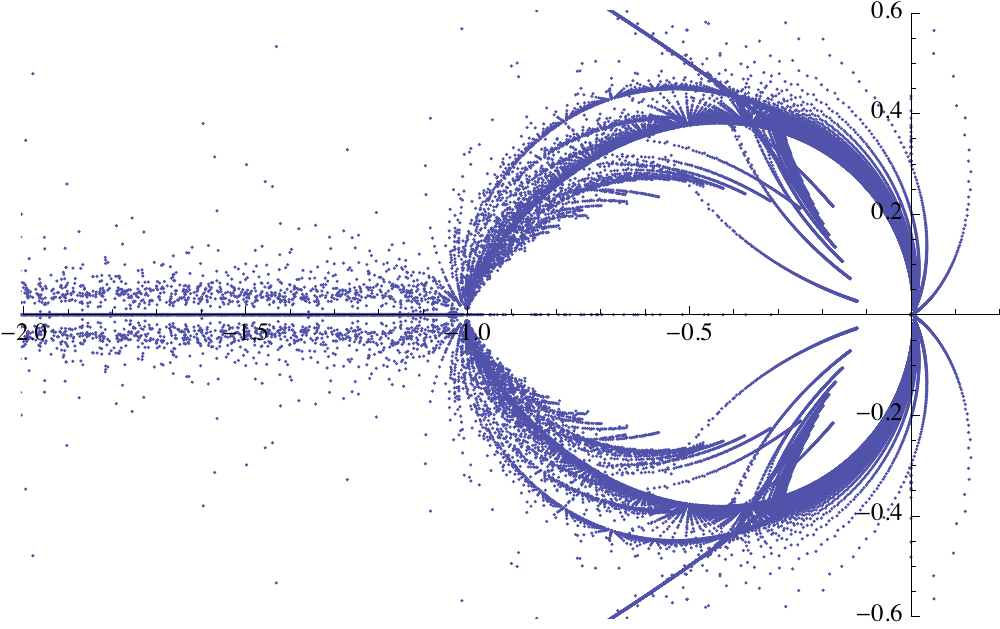,width=15em,angle=0}

\caption{$[-2,0.2] \times [-0.6,0.6]$.}
\label{fig-2-zero}
\end{minipage}
\end{figure}
}}



{{
\begin{figure}[ht]
\begin{minipage}[t]{16em}
\centering
\epsfig{file=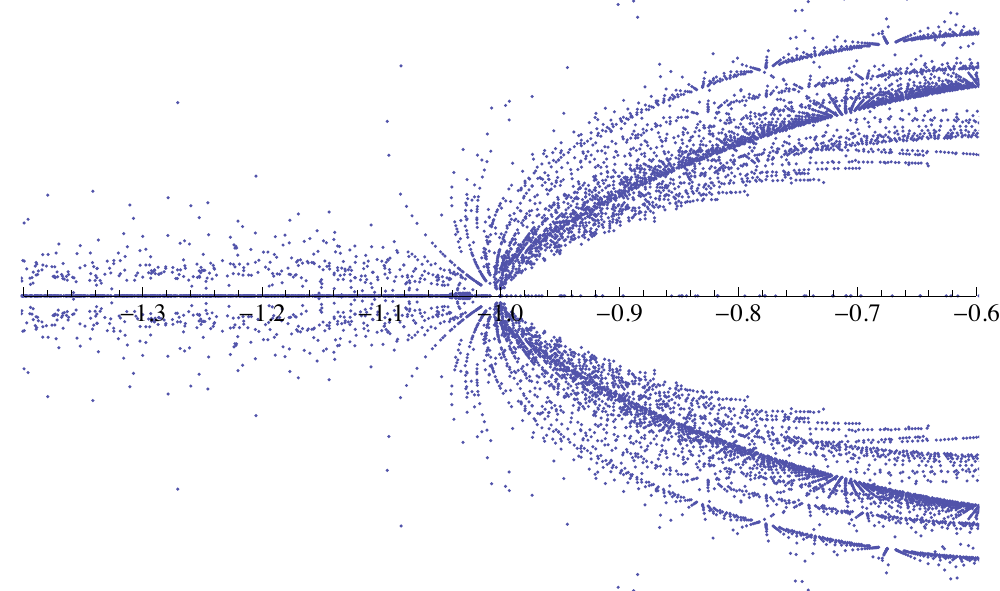,width=15em,angle=0}

\caption{$[-1.4,-0.6] \times [-0.5,0.5]$.}
\label{fig-3-zero}
\end{minipage}%
\begin{minipage}[t]{16em}
\centering
\epsfig{file=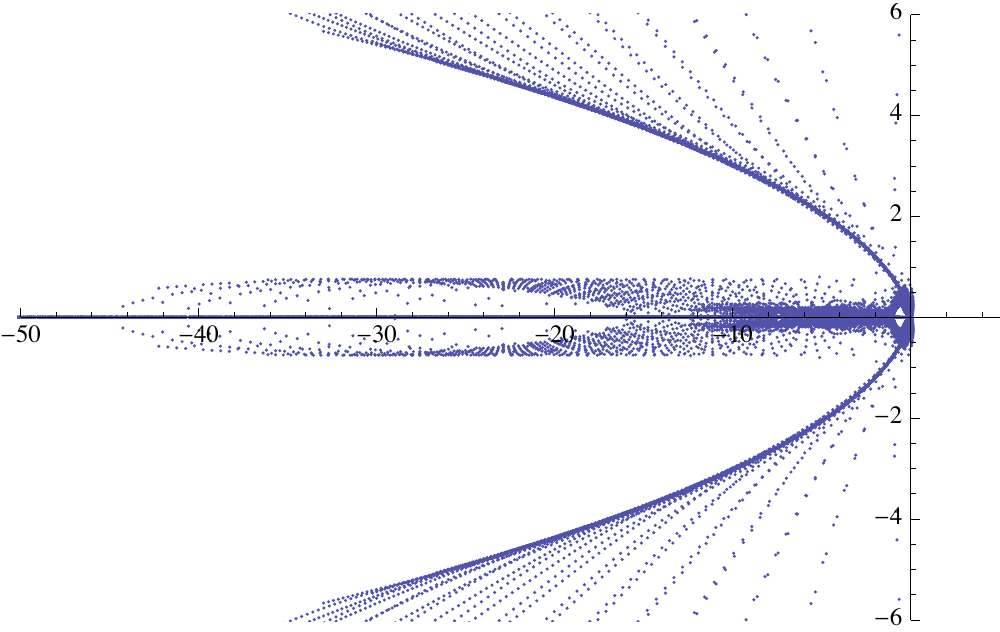,width=15em,angle=0}

\caption{$[-50,5] \times [-6,6]$.}
\label{fig-4-zero}
\end{minipage}
\end{figure}
}}

\medskip

The interesting structure depicted in figures 
\ref{fig-1-zero} to \ref{fig-4-zero} will
 be explored in future work. 

\smallskip


\noindent
{\bf Acknowledgments}. The authors wish to thank J. Silverman for 
providing information on the elliptic curve mentioned in the title and to 
A. Ayyer and A. Waldron regarding information on the quantum mechanical systems.
The authors also wish to thanks the referees for comments on an earlier 
manuscript. The third author was partially funded by
$\text{NSF-DMS } 0070567$. 

\newpage

\appendix
\section{Rows of paraffin}

The figures show all possible alkane bonds (paraffin) $C_{n}H_{2n+2}$ for 
$n = 1, \, 2, \, 3, \, 4, \, 5$.


\vspace{1cm}

\begin{center}
\begin{tikzpicture}[scale=.6, show background rectangle]

\begin{scope}
\node at (0,0) (a) {C};
\node at (0,1.5) (u1) {H};
\node at (0,-1.5) (d1) {H};
\node at (-1.5,0) (e1) {H};
\node at (1.5,0) (e2) {H};

\draw[-, thick] (a) to (u1);
\draw[-, thick] (a) to (d1);
\draw[-, thick] (a) to (e1);
\draw[-, thick] (a) to (e2);
\end{scope}

\begin{scope}[xshift=1.8in]
\node at (-1.5,0) (a1) {H};      
\node at (0,0) (a) {C};
\node at (0,1.5) (u1) {H};
\node at (0,-1.5) (u2) {H};
\node at (1.5,0) (b) {C};
\node at (1.5,1.5) (d1) {H};
\node at (1.5,-1.5) (d2) {H};
\node at (3,0) (d3) {H};

\draw[-, thick] (a) to (a1);
\draw[-, thick] (a) to (b);
\draw[-, thick] (a) to (u1);
\draw[-, thick] (a) to (u2);
\draw[-, thick] (b) to (d1);
\draw[-, thick] (b) to (d2);
\draw[-, thick] (b) to (d3);

\end{scope}

\begin{scope}[xshift=4 in]

\node at (-1.5,0) (a1) {H};      
\node at (0,0) (a) {C};
\node at (0,1.5) (u1) {H};
\node at (0,-1.5) (u2) {H};
\node at (1.5,0) (b) {C};
\node at (1.5,1.5) (d1) {H};
\node at (1.5,-1.5) (d2) {H};
\node at (3,0) (c) {C};
\node at (4.5,0) (e0) {H};
\node at (3,1.5) (e1) {H};
\node at (3,-1.5) (e2) {H};

\draw[-, thick] (a) to (a1);
\draw[-, thick] (a) to (b);
\draw[-, thick] (a) to (u1);
\draw[-, thick] (a) to (u2);
\draw[-, thick] (b) to (d1);
\draw[-, thick] (b) to (d2);
\draw[-, thick] (b) to (c);
\draw[-, thick] (c) to (e1);
\draw[-, thick] (c) to (e2);
\draw[-, thick] (c) to (e0);

\end{scope}

\draw (1,-3) node[text width=3.1cm, text justified] {$n=1:\ C H_{4}$.};
\draw (6,-3) node[text width=3.1cm, text justified] {$n=2:\ C_2 H_{6}$.};
\draw (12,-3) node[text width=3.1cm, text justified] {$n=3:\ C_3 H_{8}$.};

\end{tikzpicture}
\end{center}


\vspace{1cm}

\begin{center}
\begin{tikzpicture}[scale=.5, show background rectangle]

\draw (8,-4) node[text width=3.1cm, text justified] {$n=4:\ C_4 H_{10}$.};

\begin{scope}

\node at (0,0) (a1) {C};
\node at (1.5,0) (b) {C};
\node at (3,0) (c) {C};
\node at (4.5,0) (d) {C};
\node at (0,1.5) (u1) {H};
\node at (1.5,1.5) (u2) {H};
\node at (3,1.5) (u3) {H};
\node at (4.5,1.5) (u4) {H};
\node at (0,-1.5) (d1) {H};
\node at (1.5,-1.5) (d2) {H};
\node at (3,-1.5) (d3) {H};
\node at (4.5,-1.5) (d4) {H};
\node at (-1.5,0) (e1) {H};
\node at (6,0) (e2) {H};

\draw[-, thick] (a1) to (b);
\draw[-, thick] (c) to (b);
\draw[-, thick] (c) to (d);
\draw[-, thick] (a1) to (e1);
\draw[-, thick] (a1) to (u1);
\draw[-, thick] (a1) to (d1);
\draw[-, thick] (b) to (u2);
\draw[-, thick] (b) to (d2);
\draw[-, thick] (c) to (u3);
\draw[-, thick] (c) to (d3);
\draw[-, thick] (d) to (u4);
\draw[-, thick] (d) to (d4);
\draw[-, thick] (d) to (e2);

\end{scope}

\begin{scope}[xshift = 5 in]
      
\node at (0,0) (a) {C};
\node at (1.5,0) (b) {C};
\node at (-1.5,0) (c) {C};
\node at (0,3) (d) {C};
\node at (0,4.5) (e1) {H};
\node at (1.5,3) (e2) {H};
\node at (-1.5,3) (e3) {H};
\node at (-3,0) (f1) {H};
\node at (-1.5,-1.5) (f2) {H};
\node at (-1.5,1.5) (f3) {H};
\node at (3,0) (g1) {H};
\node at (1.5,-1.5) (g2) {H};
\node at (1.5,1.5) (g3) {H};

\draw[-, thick] (a) to (b);
\draw[-, thick] (a) to (c);
\draw[-, thick] (a) to (d);
\draw[-, thick] (d) to (e1);
\draw[-, thick] (d) to (e2);
\draw[-, thick] (d) to (e3);
\draw[-, thick] (b) to (g1);
\draw[-, thick] (b) to (g2);
\draw[-, thick] (b) to (g3);
\draw[-, thick] (c) to (f1);
\draw[-, thick] (c) to (f2);
\draw[-, thick] (c) to (f3);
\end{scope}

\end{tikzpicture}

\end{center}


\vspace{.8cm}

\begin{center}
\begin{tikzpicture}[scale=.45, show background rectangle]

\begin{scope}

\node at (0,0) (a1) {C};
\node at (1.5,0) (b) {C};
\node at (3,0) (c) {C};
\node at (4.5,0) (d) {C};
\node at (0,1.5) (u1) {H};
\node at (1.5,1.5) (u2) {H};
\node at (3,1.5) (u3) {H};
\node at (4.5,1.5) (u4) {H};
\node at (0,-1.5) (d1) {H};
\node at (1.5,-1.5) (d2) {H};
\node at (3,-1.5) (d3) {H};
\node at (4.5,-1.5) (d4) {H};
\node at (-1.5,0) (e1) {H};
\node at (6,0) (e) {C};
\node at (6,1.5) (u5) {H};
\node at (6,-1.5) (d5) {H};
\node at (7.5,0) (f5) {H};

\draw[-, thick] (a1) to (b);
\draw[-, thick] (c) to (b);
\draw[-, thick] (c) to (d);
\draw[-, thick] (a1) to (e1);
\draw[-, thick] (a1) to (u1);
\draw[-, thick] (a1) to (d1);
\draw[-, thick] (b) to (u2);
\draw[-, thick] (b) to (d2);
\draw[-, thick] (c) to (u3);
\draw[-, thick] (c) to (d3);
\draw[-, thick] (d) to (u4);
\draw[-, thick] (d) to (d4);
\draw[-, thick] (d) to (e);
\draw[-, thick] (e) to (u5);
\draw[-, thick] (e) to (d5);
\draw[-, thick] (e) to (f5);

\end{scope}

\begin{scope}[xshift = 5 in]
      
\node at (0,0) (a) {C};
\node at (1.5,0) (b) {C};
\node at (-1.5,0) (c) {C};
\node at (0,3) (d) {C};
\node at (0,4.5) (e1) {H};
\node at (1.5,3) (e2) {H};
\node at (-1.5,3) (e3) {H};
\node at (-3,0) (f1) {H};
\node at (-1.5,-1.5) (f2) {H};
\node at (-1.5,1.5) (f3) {H};
\node at (3,0) (g1) {C};
\node at (1.5,-1.5) (g2) {H};
\node at (1.5,1.5) (g3) {H};
\node at (3,1.5) (gg1) {H};
\node at (3,-1.5) (gg2) {H};
\node at (4.5,0) (gg3) {H};

\draw[-, thick] (a) to (b);
\draw[-, thick] (a) to (c);
\draw[-, thick] (a) to (d);
\draw[-, thick] (d) to (e1);
\draw[-, thick] (d) to (e2);
\draw[-, thick] (d) to (e3);
\draw[-, thick] (b) to (g1);
\draw[-, thick] (b) to (g2);
\draw[-, thick] (b) to (g3);
\draw[-, thick] (c) to (f1);
\draw[-, thick] (c) to (f2);
\draw[-, thick] (c) to (f3);
\draw[-, thick] (gg1) to (g1);
\draw[-, thick] (gg2) to (g1);
\draw[-, thick] (gg3) to (g1);
\end{scope}

\begin{scope}[xshift = 9 in]
      
\node at (0,0) (a) {C};
\node at (1.5,0) (b) {C};
\node at (-1.5,0) (c) {C};
\node at (0,3) (d) {C};
\node at (0,-3) (e) {C};
\node at (0,4.5) (e1) {H};
\node at (1.5,3) (e2) {H};
\node at (-1.5,3) (e3) {H};
\node at (-3,0) (f1) {H};
\node at (-1.5,-1.5) (f2) {H};
\node at (-1.5,1.5) (f3) {H};
\node at (3,0) (g1) {H};
\node at (1.5,-1.5) (g2) {H};
\node at (1.5,1.5) (g3) {H};
\node at (-1.5,-3) (ee1) {H};
\node at (1.5,-3) (ee2) {H};
\node at (0,-4.5) (ee3) {H};

\draw[-, thick] (a) to (b);
\draw[-, thick] (a) to (c);
\draw[-, thick] (a) to (d);
\draw[-, thick] (d) to (e1);
\draw[-, thick] (d) to (e2);
\draw[-, thick] (d) to (e3);
\draw[-, thick] (b) to (g1);
\draw[-, thick] (b) to (g2);
\draw[-, thick] (b) to (g3);
\draw[-, thick] (c) to (f1);
\draw[-, thick] (c) to (f2);
\draw[-, thick] (c) to (f3);

\draw[-, thick] (a) to (e);
\draw[-, thick] (e) to (ee1);
\draw[-, thick] (e) to (ee2);
\draw[-, thick] (e) to (ee3);

\end{scope}

\centering

\draw (15,-6) node[text width=3.1cm, text justified] {$n=5:\ C_5 H_{12}$.};

\end{tikzpicture}

\end{center}




\end{document}